\documentclass[a4paper,12pt,leqno]{article}

\usepackage[latin1]{inputenc} 
\usepackage[T1]{fontenc}      
\usepackage{geometry}         
\usepackage[francais, english]{babel}  
\usepackage{amsmath}
\usepackage{amsfonts}
\usepackage{amssymb}
\usepackage{fancyhdr}
\usepackage{url}
\usepackage{mathrsfs}
\usepackage{euscript}
\usepackage{enumerate}
\usepackage{cite}
\usepackage{textcomp}
\usepackage[all]{xy}
\usepackage{mathtools}

\newtheorem{df}{Definition}[section]
\newtheorem{nt}[df]{Notation}

\newtheorem{prop}[df]{Proposition}
\newtheorem{thm}[df]{Theorem}
\newtheorem{dfn}{Definition}
\newtheorem{thrm}[dfn]{Theorem}
\newtheorem{conj}[dfn]{Conjecture}
\newtheorem{lem}[df]{Lemma}

\newcommand{\prf}{\noindent \textit{Proof}}
\newtheorem{rmk}[df]{Remark}

\newcommand{\R}{\mathbb{R}}

\newcommand{\XD}{X\backslash D}

\newcommand{\vareps}{\varepsilon}

\newcommand{\vl}{\operatorname{Vol}}

\newcommand{\scal}{\mathbf{s}}
\newcommand{\id}{\operatorname{id}}

\newcommand{\cqfd}{ \hfill $\square$ }

\newcommand{\sbar}{\overline{\scal}}

\newcommand{\varparallel}{\mathrel{\!/\mkern-5mu/\!}}
\newcommand{\hD}{\mathfrak{h}_{\varparallel}^D}

\title{\large NOTE ON POINCARÉ TYPE KÄHLER METRICS AND FUTAKI CHARACTERS} 

\vspace{5pt}
\author{\normalsize\textsc{Hugues AUVRAY}}
\date{}

\begin{document}

\makeatletter
\renewcommand%
   {\section}%
   {%
   \@startsection{section}%
      {1}
      {0mm}%
      {\baselineskip}
      {0.5\baselineskip}%
      {\sc\large\centering}
   }%
\makeatother

\makeatletter
\renewcommand%
   {\subsubsection}%
   {%
   \@startsection{subsubsection}%
      {1}
      {0mm}%
      {1.25\baselineskip}
      {0.25\baselineskip}%
      {\sf\normalsize}
   }%
\makeatother

\makeatletter
\renewcommand%
   {\paragraph}%
   {%
   \@startsection{paragraph}%
      {4}%
      {0mm}%
      {0mm}%
      {-5pt}%
      {\it\normalsize}%
   }%
\makeatother

 \maketitle
 
 \begin{abstract}
  A Poincaré type Kähler metric on the complement $\XD$ of a simple normal crossing divisor $D$,  
  in a compact Kähler manifold $X$, is a Kähler metric on $\XD$ with cusp singularity along $D$. 
  We relate the \textit{Futaki character for holomorphic vector fields parallel to the divisor}, 
  defined for any fixed Poincaré type Kähler class, 
  to the classical Futaki character for the relative smooth class. 
  As an application we express a numerical obstruction to the existence of extremal Poincaré type Kähler metrics, 
  in terms of mean scalar curvatures and Futaki characters. 
 \end{abstract}

 \begin{center}
   \rule{3cm}{0.25pt}
 \end{center}

~

\section*{Introduction}
 A basic fact in Kähler geometry is the independence of the de Rham class of the Ricci form from the background metric 
 on a compact Kähler manifold: 
 it is always $-2\pi c_1(K)$, with $c_1(K)$ the first Chern class of the canonical line bundle.  
 This topological invariance constitutes the first obstacle for a compact Kähler manifold to admit a Kähler-Einstein metric: 
 the Chern class in question must then have a sign, 
 which, if definite, forces Kähler-Einstein metrics to lie 
 in a consequently fixed Kähler class. 
 
 When $c_1(K)>0$, a (unique) Kähler-Einstein metric was obtained by Aubin and Yau, 
 and Bochner's technique then rules out the existence of non-trivial holomorphic vector fields .
 Conversely, in the opposite case $c_1(K)<0$, the so-called \textit{``Fano case''}, 
 non-trivial holomorphic vector fields may exist, 
 and the existence of a Kähler-Einstein metric, which does not always hold, is noticeably more involved. 
 More precisely, 
 in this case -- and, respectively, on any compact Kähler manifold -- 
 non-trivial holomorphic vector fields bring a constraint to the existence of a Kähler-Einstein metric 
 -- respectively, of a constant scalar curvature metric Kähler metric in a fixed Kähler class.  
 If indeed such a canonical metric exists, a numerical function, the \textit{Futaki character} \cite{fut}, 
 defined on the Lie algebra of holomorphic vector fields 
 -- and depending only on the Kähler class under study --, 
 has to vanish identically . 

 The Futaki character was later generalised by Donaldson 
 to polarised manifolds, 
 into an numerical function defined on \textit{test-configurations}, 
 which generalise the concept of (the action of) holomorphic vector fields \cite{donx}. 
 In the lines of suggestions by Yau \cite{yau}, 
 and after Tian's \textit{special degenerations} \cite{tia}, 
 test-configurations and their Donaldson-Futaki invariants are meant to reveal the link between 
 algebro-geometric stability of the manifold, and existence of a Kähler-Einstein/constant scalar curvature Kähler metric: 
  \begin{conj}[Yau-Tian-Donaldson]
   A polarised manifold $(X,L)$ admits a constant scalar Kähler curvature metric in $2\pi c_1(L)$ 
   if, and only if, 
   $(X,L)$ is ``K-stable'', that is: 
   the Donaldson-Futaki invariant is nonpositive (negative) for any (non-trivial) test-configuration. 
  \end{conj}
 The ``only if'' direction is now established \cite{mab, sto}; 
 the ``if'' direction is still a very active area of research, 
 and has recently been solved for Kähler-Einstein metrics in the Fano case, 
 i.e. when $L=-K_X$ is ample, 
 see \cite{cds1,cds2,cds3} and \cite{tia2}. 
 
 In a related scope, the aim of this note is, 
 after restriction to the relevant set of holomorphic vector fields, 
 to generalise the Futaki character to a certain class of singular metrics 
 on a compact manifold. 
 Namely, fixing a simple normal crossing divisor $D$ in a compact Kähler manifold $(X,J, \omega_X)$, 
 we recall the definition \textit{Poincaré type Kähler metrics on $\XD$}, 
 following \cite{tian-yau1, wu2, auv1}: 
  \begin{dfn}
   A smooth positive $(1,1)$-form $\omega$ on $\XD$ is called a 
   \textbf{Poincaré type Kähler metric on $\XD$} 
   if: 
   on every open subset $U$ of coordinates $(z^1,\dots, z^m)$ in $X$, 
   in which $D$ is given by $\{z^1\cdots z^j=0\}$, 
   $\omega$ is mutually bounded with 
    \begin{equation*}
     \omega_{U}^{\rm mdl}
     :=\frac{i dz^1\wedge d\overline{z^1}}{|z^1|^2\log^2(|z^1|^2)}
     +\cdots+\frac{i dz^j\wedge d\overline{z^j}}{|z^j|^2\log^2(|z^j|^2)}
     +i dz^{j+1}\wedge d\overline{z^{j+1}}
     +\cdots + i dz^{m}\wedge d\overline{z^{m}}, 
    \end{equation*}
   and has bounded derivatives at any order for this model metric. 
   
   We say moreover that \textbf{$\omega$ is of class $[\omega_X]$} 
   if $\omega=\omega_X+dd^c\varphi$ for some $\varphi$ smooth on $\XD$, 
   with $\varphi=\mathcal{O}\big(\sum_{\ell=1}^j \log[-\log(|z^{\ell}|^2)]\big)$ in the above coordinates 
   and $d\varphi$ bounded at any order for $\omega^{\rm mdl}_U$. 
   We then set: $\omega \in \mathscr{M}_{[\omega_X]}^D$. 
  \end{dfn}

  Metrics of $\mathscr{M}_{[\omega_X]}^D$ are complete, 
  with finite volume (equal to that of $X$ for smooth Kähler metrics of class $[\omega_X]$); 
  they also share a common mean scalar curvature, 
  which differs from that attached to smooth Kähler metrics of class $[\omega_X]$. 
  
  Restricting our attention to the set $\hD$ of holomorphic vector fields with their normal component vanishing along $D$ 
  -- these are the holomorphic vector fields bounded (at any order) for \textit{any} Poincaré type Kähler metric on $\XD$ --, 
  we define in a way similar to the compact case a \textit{Poincaré type Futaki character} $\mathscr{F}^{D}_{[\omega_X]}$, 
  computed with metrics of $\mathscr{M}_{[\omega_X]}^D$, and depending only on this Poincaré class. 
  
 ~
  
  \noindent 
  \textit{Results. --- }
  This Poincaré type Futaki character differs generally from the restriction of the usual smooth Futaki character 
  $\mathscr{F}_{[\omega_X]}$ to $\hD$. 
  More precisely, our first main result is a formula giving the precise relation between these two invariants: 
   \begin{thrm}  \label{thm_mainthm1}
    Let $\mathsf{Z}\in \hD$, with Riemannian gradient potential $f$ for $\omega_X(\cdot,J\cdot)$. 
    Then 
     $\mathscr{F}^{D}_{[\omega_X]}(\mathsf{Z})
      =\mathscr{F}_{[\omega_X]}(\mathsf{Z})+\sum_{j=1}^N\int_{D_j} f\,\frac{(\omega_X|_{D_j})^{m-1}}{(m-1)!}$. 
   \end{thrm}
  The gradient potential refers to that of the Hodge decomposition of (the dual 1-form) of $\mathsf{Z}$. 
  
  As an application of Poincaré type Futaki character, 
  in the framework of finding necessary conditions on canonical Kähler metrics, 
  we provide the following \textit{numerical constraint on the existence of extremal metrics of Poincaré type on $\XD$}: 
   \begin{thrm} \label{thm_mainthm2}
    Assume that there exists an extremal metric in $\mathscr{M}^{D}_{[\omega_X]}$, 
    and denote by $\mathsf{K}$ the Riemannian gradient of its scalar curvature. 
    Then for all $j\in \{1, \dots, N\}$ indexing an irreducible component of $D$, 
    one has: 
     \begin{equation}   \label{eqn_constraint}
      \sbar^D< \sbar_{D_j}^{D^j}
               +\frac{1}{4\pi\vl(D_j)}\big(\mathscr{F}^{D-D_j}_{[\omega_X]}(\mathsf{K})-\mathscr{F}^{D}_{[\omega_X]}(\mathsf{K})\big),  
     \end{equation}
    where $\mathscr{F}^{D-D_j}_{[\omega_X]}(K)$ is the Futaki character for Poincaré type metrics 
    on $X\backslash (D-D_j)=X\backslash \sum_{\ell\neq j}D_{\ell}$ of class $[\omega_X]$, 
    $\sbar^D$ is the mean scalar curvature attached to $\mathscr{F}^{D}_{[\omega_X]}(K)$, 
    and $\sbar_{D_j}^{D^j}$ that attached to $\mathscr{F}^{D^j}_{[\omega_X]|_{D_j}}(K)$, 
    the class of Poincaré type metrics on $D_j\backslash (D-D_j)|_{D_j}$ of class $[\omega_X]|_{D_j}$. 
   \end{thrm}
   
   An extremal metrics is a Kähler metrics such that the Riemannian gradient of its scalar curvature is holomorphic. 
   Constraint \eqref{eqn_constraint} is a reformulation of that of \cite[Prop.4.5]{auv3}, 
   and as such, extends that on the existence of constant scalar curvature metrics in $\mathscr{M}^{D}_{[\omega_X]}$ 
   of \cite{auv2}, which states as: 
   $\sbar^D< \sbar_{D_j}^{D^j}$ for all $j=1,\dots, N$. 
   As in the compact case moreover, 
   by construction and invariance on $\mathscr{M}^{D}_{[\omega_X]}$, 
   $\mathscr{F}^{D}_{[\omega_X]}$ vanishes identically if there exists a constant scalar curvature metric in $\mathscr{M}^{D}_{[\omega_X]}$, 
   and conversely, its vanishing forces possible extremal metrics of $\mathscr{M}^{D}_{[\omega_X]}$ 
   to have constant scalar curvature. 
   By contrast, 
   the interest of Theorem \ref{thm_mainthm2} is to provide constraints on extremal metrics independently of such a vanishing.  
   
   Finally, Donaldson-Futaki invariants are already considered in \cite{cds0} 
   which take into account the contribution of a divisor. 
   These are used in the context of \textit{Kähler metrics with conical singularities on polarised manifolds}, 
   and the divisor term of the invariant comes with coefficient $(1-\beta)$, 
   with $2\pi\beta$ the angle of the cone singularity. 
   In view of Theorem \ref{thm_mainthm1}, 
   the Poincaré type Futaki invariant might thus be viewed 
   -- at the level of holomorphic vector fields rather than at that of test-configurations -- 
   as the limit when the conical singularity angle goes to $0$, 
   that is, roughly speaking, when \textit{cones} become \textit{cusps}. 
   
 ~
 
  \noindent
  \textit{Organisation of the article. --- }
   This note is divided into three parts. 
   In the first part, 
   we analyse holomorphic vector fields parallel to a divisor, 
   with respect to Poincaré type Kähler metrics on the complement of this divisor. 
   We see in particular that a Hodge decomposition analogous to that of the compact case still holds for such metrics and such vector fields. 
   This allows us in Section \ref{subsec_PFchar} to define the Poincaré type Futaki character, 
   as an \textit{invariant} of a given Poincaré type Kähler class. 
   
   Theorem \ref{thm_mainthm1}, under a slightly more general version, 
   is stated in Section \ref{subsec_sttmt} (Proposition \ref{prop_futakis}); 
   it is proven in Section \ref{subsec_PrfPropFutakis}, 
   and the final section \ref{subsec_PrfKeyLm} of Part 2 is devoted to a key technical lemma (Lemma \ref{lem_keylm}) 
   used in Section \ref{subsec_PrfPropFutakis}. 
   
   In Part 3 we state and prove Theorem \ref{thm_mainthm2}: 
   a useful extension (Proposition \ref{prop_futakis2}) of \ref{thm_mainthm1} 
   (to asymptotically product Poincaré type metrics when the divisor is smooth) 
   is given in Section \ref{subsec_ExtPropFtks};  
   Theorem \ref{thm_mainthm2} is then proven in Section \ref{subsec_ExtCnstrnt} (Theorem \ref{thm_ExtCnstt}), 
   first in the smooth divisor case using Proposition \ref{prop_futakis2}, 
   then in the simple normal crossing case. 
   Notice that both steps use the asymptotic properties of extremal Poincaré type metrics obtained in \cite{auv3}.

 \subsubsection*{Acknowledgements} 
  I am very thankful to Vestislav Apostolov and Yann Rollin for the informal discussions I had with them, 
  which led me to the material and results presented here.  
  
 ~

 \textbf{\textit{In all this note, $X$ is a compact Kähler manifold, 
 and $D\subset X$ a simple normal crossing divisor, 
 the decomposition into irreducible smooth components of which we write as $\sum_{j=1}^N D_j$.} }

\section{The Futaki character of a Poincaré class}     \label{sctn_PFchar}
 \subsection{Hodge decomposition of vector fields parallel to the divisor}   \label{subsec_HdgDec}
 \noindent
 \textit{Reminder: the compact case. ---}
  Fix a smooth Kähler form $\omega_X$ on $X$, 
  of associated Riemannian metric $g_X$. 
  Given any real holomorphic vector field $\mathsf{Z}$ -- ``$\mathsf{Z}\in\mathfrak{h}$'' --, 
  it is well-known that its $g_X$-dual 1-form $\xi^{\mathsf{Z}}$, that is, $\mathsf{Z}^{\sharp_{g_X}}$, 
  enjoys the following decomposition:
   \begin{equation}   \label{eqn_Hodgedec}
    \xi^{\mathsf{Z}} = \xi^{\mathsf{Z}}_{\rm harm} + df^{\mathsf{Z}}_{\omega_X} + d^ch^{\mathsf{Z}}_{\omega_X}, 
   \end{equation}
  into harmonic, $d$- and $d^c$-parts;  
  these are uniquely determined, 
  provided that $f^{\mathsf{Z}}_{\omega_X}$ and $h^{\mathsf{Z}}_{\omega_X}$ are taken with null mean against $\omega_X^m$. 
  
  Decomposition \eqref{eqn_Hodgedec} is called the \textit{(dual) Hodge decomposition} of $\mathsf{Z}$.  
  Given moreover any other smooth metric $\tilde{\omega}=\omega_X+dd^c\varphi$ of $\mathscr{M}_{[\omega_X]}$, 
  and setting $\tilde{\xi}^{\mathsf{Z}}$ for the dual 1-form of $\mathsf{Z}$ with respect to $\tilde{\omega}(\cdot,J\cdot)$, 
  its Hodge decomposition is: 
   \begin{equation*}   
    \tilde{\xi}^{\mathsf{Z}} = \xi^{\mathsf{Z}}_{\rm harm} + d\big(f^{\mathsf{Z}}_{\omega_X}+\mathsf{Z}\cdot\varphi\big) 
                                                           + d^c\big(h^{\mathsf{Z}}_{\omega_X}-(J\mathsf{Z})\cdot\varphi\big) 
   \end{equation*}
  see \cite[Lemma 4.5.1]{gau}-- notice in particular that the harmonic part remains unchanged at the level of 1-forms; 
  recall that on compact Kähler manifolds, the space of harmonic 1-forms is \textit{independent of the Kähler metric}.

 ~
 
  \noindent 
  \textit{Extension to Poincaré type Kähler metrics for vector fields parallel to a divisor. --- }
   Consider now the simple normal crossing divisor $D=\sum_{j=1}^N D_j$ in $X$. 
   The normal crossing assumption can be expressed as follows: 
   given any $p\in (D_{j_1}\cap\cdots\cap D_{j_k})\backslash (D_{\ell_1}\cup\cdots\cup D_{\ell_{N-k}})$, 
   with $\{j_1,\dots,j_k\}\sqcup\{\ell_1\dots,\ell_{N-k}\}=\{1,\dots,N\}$, 
   one can find in $X$ an open neighbourhood $U$ of $p$, 
   of holomorphic coordinates $(z^1,\dots,z^m)$ 
   such that $U\cap D_{j_s}=\{z^s=0\}$ for $s=1,\dots,k$ (in particular, $k\leq m$). 
   
   We define a restricted class of holomorphic vector fields on $X$, 
   the use of which is natural when working with Poincaré type Kähler metrics on $\XD$:
  \begin{df}
   Let $\mathsf{Z}\in \mathfrak{h}$. 
   We say that $\mathsf{Z}$ \textbf{is parallel to $D$}, denoted $\mathsf{Z}\in \hD$, 
   if: 
   writing $\mathsf{Z}$ as $\mathfrak{Re}\big(f_1\frac{\partial}{\partial z^1}+\cdots +f_m\frac{\partial}{\partial z^m} \big)$ 
   in local holomorphic coordinates as above, 
   one has $f_1\equiv\cdots\equiv f_s\equiv0$ on $D$. 
   
   Given $j\in \{1,\dots, N\}$, 
   we then define \textbf{the restriction $\mathsf{Z}|_{D_j}$ of $\mathsf{Z}$ to $D_j$} 
   by setting locally 
   $\mathsf{Z}_{D_j} = \mathfrak{Re}\big(f_2|_{D_j}\frac{\partial}{\partial z^2}+\cdots +f_m|_{D_j}\frac{\partial}{\partial z^m}\big)$, 
   whenever $j_1=j$ in the above coordinates. 
  \end{df}

  One checks in particular that the definition of $\mathsf{Z}_{D_j}$ is independent of the choice of 
  holomorphic coordinates, as long as the first coordinate is a local equation of $D_j$; 
  one also checks easily that $\hD$ is a Lie subalgebra of $\mathfrak{h}$. 
  
  Holomorphic vector fields parallel to $D$ are relevant when working with Poincaré type Kähler metrics on $D$ 
  for the following reason (see the proof of Lemma 5.2 in \cite{auv1}): 
  \textit{any holomorphic vector field on $\XD$ which is bounded -- or actually, $L^2$ -- with respect to 
  a Poincaré type Kähler metric on $\XD$ extends to a holomorphic vector field on $X$, parallel to $D$}. 
  Conversely, 
  any holomorphic vector field on $X$ parallel to $D$ gives on $\XD$ a vector field bounded at any order for any Poincaré type metric 
  on $\XD$.

  ~

  We now provide \textit{a Hodge decomposition of holomorphic vector fields parallel to $D$ 
  with respect to Poincaré type Kähler metrics on $\XD$}, 
  analogous to the decomposition of the compact setting: 
   \begin{prop}  \label{prop_HodgeDec}
    Let $\mathsf{Z}\in \hD$, and let $\omega=\omega_X+dd^c\varphi \in \mathscr{M}_{[\omega_X]}^{D}$. 
    Let $\xi^{\mathsf{Z}}_{\varphi}$ be the dual 1-form of $\mathsf{Z}$ with respect to $\omega(\cdot,J\cdot)$. 
    Then 
     \begin{equation}   \label{eqn_Hodgedecphi}
      \xi^{\mathsf{Z}}_{\varphi} = \xi^{\mathsf{Z}}_{\rm harm} + d\big(f^{\mathsf{Z}}_{\omega_X}+\mathsf{Z}\cdot\varphi\big) 
                                                           + d^c\big(h^{\mathsf{Z}}_{\omega_X}-(J\mathsf{Z})\cdot\varphi\big) 
     \end{equation}
    on $\XD$,  
    with the same harmonic part $\xi^{\mathsf{Z}}_{\rm harm}$ as in the compact case, 
    and this decomposition is unique. 
    Moreover, 
     \begin{equation*}
      \int_{\XD} \big(f^{\mathsf{Z}}_{\omega_X}+\mathsf{Z}\cdot\varphi\big) \omega^m =
      \int_{\XD} \big(h^{\mathsf{Z}}_{\omega_X}-(J\mathsf{Z})\cdot\varphi\big) \omega^m = 0. 
     \end{equation*}
   \end{prop}
   
  The uniqueness we state here is understood as follows: 
  if $\xi^{\mathsf{Z}}_{\varphi}= \alpha + d\beta + d^c\gamma$ 
  with $\alpha$ harmonic on $\XD$, 
  and $\alpha$, $\beta$, $\gamma$ bounded for $\omega$ of Poincaré type, 
  then $\alpha=\xi^{\mathsf{Z}}_{\rm harm}$, 
  and 
  $\beta=f^{\mathsf{Z}}_{\omega_X}+\mathsf{Z}\cdot\varphi$ and $\gamma=h^{\mathsf{Z}}_{\omega_X}-(J\mathsf{Z})\cdot\varphi$ 
  up to a constant. 
  This justifies: 
   \begin{nt} 
    With the notations of Proposition \ref{prop_HodgeDec},
    we set 
     \begin{equation*}
      f^{\mathsf{Z}}_{\omega} = f^{\mathsf{Z}}_{\omega_X}+\mathsf{Z}\cdot\varphi
       \quad\text{ and }\quad
      h^{\mathsf{Z}}_{\omega} = h^{\mathsf{Z}}_{\omega_X}-(J\mathsf{Z})\cdot\varphi. 
     \end{equation*}

   \end{nt}

  ~
  
 \noindent
 \prf \textit{of Proposition \ref{prop_HodgeDec}. --- }
  With the notations of the statement, 
  we first prove that equality \eqref{eqn_Hodgedecphi} holds on $\XD$.  
  This identity is purely local; 
  it is thus sufficient to establish it for any Kähler metric equal to $\omega$ in the neighbourhood 
  of any given point of $\XD$. 
  More concretely, 
  as $\omega=\omega_X+dd^c\varphi$ is of Poincaré type, 
  local analysis provides that $\varphi\to -\infty$ near $D$. 
  Consider a convex function $\chi:\R\to\R$, 
  with $\chi(t)=0$ if $t\leq -1$, $\chi(t)=t$ if $t\geq 1$ -- and thus $0\leq \chi'(t) \leq 1$ for all $t$. 
  Given $K\in \R$, one now easily checks that 
   \begin{equation*}   \label{eqn_omegaK}
    \omega_K := \omega_X +dd^c\big(\chi \circ (\varphi + K)\big)
   \end{equation*}
  is a smooth metric on $X$, 
  equal to $\omega$ on $\{\varphi\geq 1-K\}$ (compact in $\XD$), 
  and to $\omega_X$ on $\{\varphi\leq -(K+1)\}$. 
  Now \eqref{eqn_Hodgedecphi} follows on $\{\varphi > 1-K\}$ by the smooth case of Hodge decomposition applied to $\omega_K$, 
  thus on all $\XD$ by letting $K\to\infty$. 
  
  Observe that $\xi^{\mathsf{Z}}_{\rm harm}$ is still harmonic with respect to $\omega$; 
  again, this condition is local, 
  implied, thanks to the Kähler identities, 
  by the closedness and the $d^c$-closedness of $\xi^{\mathsf{Z}}_{\rm harm}$. 
  These latter conditions are independent of the Kähler metric, 
  and indeed implied by the harmonicity of $\xi^{\mathsf{Z}}_{\rm harm}$ for the smooth $\omega_X$, as $X$ is compact. 
  
  As $\xi^{\mathsf{Z}}_{\rm harm}$ is bounded for $\omega_X$, 
  it is so for $\omega$, which dominates $\omega_X$. 
  Similarly, $f^{\mathsf{Z}}_{\omega_X}$ and $h^{\mathsf{Z}}_{\omega_X}$ are bounded at any order for $\omega_X$ hence for $\omega$, 
  and as $\mathsf{Z}$ is parallel to $D$, 
  it is bounded at any order for $\omega$, as well as $d\varphi$ by definition; 
  consequently, $f^{\mathsf{Z}}_{\omega}$ 
  and $h^{\mathsf{Z}}_{\omega}$ are bounded at any order for $\omega$. 
  From this the uniqueness of Hodge decomposition easily follows.  
  Write $\xi^{\mathsf{Z}}_{\varphi}=\alpha+d\beta+d^c\gamma$ with $\alpha$, $\beta$, $\gamma$ as above.  
  As $d\alpha=d^c\alpha=0$ 
  ($\alpha$ is bounded and harmonic for $\omega$ hence bounded at any order by uniform ellipticity in quasi-coordinates, 
  and one can thus integrate by parts without boundary terms), 
  one gets $dd^c(f^{\mathsf{Z}}_{\omega}-\beta)
             =dd^c (h^{\mathsf{Z}}_{\omega}-\gamma)=0$ on $\XD$. 
  Therefore $f^{\mathsf{Z}}_{\omega}-\beta$ and $h^{\mathsf{Z}}_{\omega}-\gamma$ 
  are constant (use e.g. Yau's maximum principle \cite[p.406]{wu2}) as wanted, 
  and thus $\alpha=\xi^{\mathsf{Z}}_{\rm harm}$. 
  
  We are left with the mean assertion on $f^{\mathsf{Z}}_{\omega}$ 
  and $h^{\mathsf{Z}}_{\omega}$. 
  For $t\in[0,1]$, set $\omega_t=\omega_X+tdd^c\varphi$, $f_t=f^{\mathsf{Z}}_{\omega_X}+t(\mathsf{Z}\cdot\varphi)$, 
  and consider the function $t\mapsto \int_{\XD} f_t\,\omega_t^m$. 
  Thanks to the growths near $D$, 
  this function is smooth, 
  with derivative $\int_{\XD} (\mathsf{Z}\cdot\varphi)\,\omega_t^m + m\int_{\XD} f_t dd^c\varphi\wedge\omega_t^{m-1}$.  
  Now,  
   \begin{equation*}
    m\int_{\XD} f_t dd^c\varphi\wedge\omega_t^{m-1}=
    -m\int_{\XD} df_t\wedge d^c\varphi\wedge\omega_t^{m-1}=
    -\int_{\XD} \langle df_t, d\varphi\rangle_{\omega_t}\, \omega_t^{m} 
   \end{equation*}
  (no boundary terms). 
  On the other hand, we now know that for all $t$, 
  $\xi^{\mathsf{Z}}_{t\varphi}=\xi^{\mathsf{Z}}_{\rm harm} + d\big(f^{\mathsf{Z}}_{\omega_X}+t(\mathsf{Z}\cdot\varphi)\big) 
                                                + d^c\big(h^{\mathsf{Z}}_{\omega_X}-t(J\mathsf{Z})\cdot\varphi\big)$.  
  Notice that  
  $\int_{\XD} \langle \xi^{\mathsf{Z}}_{\rm harm}, d\varphi\rangle_{\omega_t}\, \omega_t^{m}
   = \int_{\XD} \varphi(\delta_{\omega_t}\xi^{\mathsf{Z}}_{\rm harm})\, \omega_t^{m}=0$ 
  ($\xi^{\mathsf{Z}}_{\rm harm}$ is co-closed for $\omega_t$, 
  as $\delta_{\omega_t}=\Lambda_{\omega_t}d^c$ on 1-forms), 
  and $\int_{\XD} \langle d^c\big(h^{\mathsf{Z}}_{\omega_X}-t(J\mathsf{Z})\cdot\varphi\big), d\varphi\rangle_{\omega_t}\, \omega_t^{m}=
       -m\int_{\XD} d\big(h^{\mathsf{Z}}_{\omega_X}-t(J\mathsf{Z})\cdot\varphi\big)\wedge d\varphi \wedge \omega_t^{m-1}=0$ 
  (the integrand is closed, and there are no boundary terms). 
  This way 
   $\int_{\XD} \langle df_t, d\varphi\rangle_{\omega_t}\, \omega_t^{m}
    =\int_{\XD} \langle \xi^{\mathsf{Z}}_{t\varphi}, d\varphi\rangle_{\omega_t}\, \omega_t^{m}
    =\int_{\XD} (\mathsf{Z}\cdot \varphi) \, \omega_t^m$, hence: 
  $\int_{\XD} f_t\,\omega_t^m$ is constant, 
  which gives (take $t=0,1$): 
  $\int_{\XD} f^{\mathsf{Z}}_{\omega}\, \omega^m =
   \int_{\XD} f^{\mathsf{Z}}_{\omega_X}\,\omega_X^m =0$. 
  The mean of $h^{\mathsf{Z}}_{\omega}$ against $\omega^m$ 
  is seen to vanish likewise. 
 \cqfd

 \subsection{The Poincaré type Futaki character}   \label{subsec_PFchar}
  \noindent 
  \textit{Definition. --- }
  We can now generalise to Poincaré type Kähler metrics/classes, 
  and holomorphic vector fields parallel to the divisor, 
  a well-known invariant \cite{fut} of compact Kähler manifolds: 
   \begin{df}
    For $\mathsf{Z}\in \hD$ and $\omega\in \mathscr{M}_{[\omega_X]}^D$, 
    we call \textbf{Poincaré type Futaki character of $\mathsf{Z}$ with respect to $D$} the quantity
     \begin{equation}   \label{eqn_dfFutaki}
      \mathscr{F}^D_{[\omega_X]}(\mathsf{Z}) = \int_{\XD} \scal(\omega) f^{\mathsf{Z}}_{\omega}\, \frac{\omega^m}{m!}.
     \end{equation}
   \end{df}
  Here, $\scal(\omega)$ denotes the \textit{(Riemannian) scalar curvature} of $\omega$, 
  that one can compute for instance via: $\scal(\omega)\frac{\omega^m}{m!}=2\varrho(\omega)\wedge\frac{\omega^{m-1}}{(m-1)!}$, 
  with $\varrho(\omega)$ the Ricci-form of $\omega$. 
  
  ~
  
  \noindent 
  \textit{Independence from the reference metric. --- }  
  As terminology and notation suggest, 
  this Poincaré type Futaki character does not depend on $\omega$ of class $[\omega_X]$, 
  \textit{provided it is of Poincaré type}: 
   \begin{prop}   \label{prop_InvceFutaki}
    Let $\tilde{\omega}$ be any Poincaré type metric in $\mathscr{M}_{[\omega_X]}^D$, 
    and $\mathsf{Z}\in \hD$. 
    Then $\mathscr{F}^D_{[\omega_X]}(\mathsf{Z}) 
           = \int_{\XD} \scal(\tilde{\omega}) f^{\mathsf{Z}}_{\tilde{\omega}}\, \frac{\tilde{\omega}^m}{m!}$. 
   \end{prop}
   
  ~
 
  Observe nonetheless that we take $\tilde{\omega}$ of Poincaré type in this proposition; 
  the relation between the usual smooth Futaki character, 
  and our Poincaré type Futaki character, is the purpose of next part. 
  For now, let us address the proof. 
  
  ~
  
  \noindent
  \prf \textit{of Proposition \ref{prop_InvceFutaki}. --- } 
   Take $\mathsf{Z}\in \hD$. 
   Fix $\omega=\omega_X+dd^c\varphi$ and $\tilde{\omega}=\omega_X+dd^c\tilde{\varphi}$ in $\mathscr{M}^D_{[\omega_X]}$, 
   and for $t\in [0,1]$, 
   set $\omega_t=(1-t)\omega+t\tilde{\omega}=\omega_X+dd^c\varphi_t$, 
   $\varphi_t=(1-t)\varphi+t\tilde{\varphi}$; 
   the $\omega_t$ are metrics of Poincaré type, uniformly bounded below by $c\omega$, say. 
   As a consequence, 
   the $\scal(\omega_t)$ are uniformly bounded, at any order for $\omega$, 
   and for all $t_0\in[0,1]$, 
   $\scal(\omega_{t})=\scal(\omega_{t_0})+(t-t_0)\dot{\scal}_{t_0}+(t-t_0)^2w_{t_0,t}$, 
   with $\dot{\scal}_{t_0}=-\frac{1}{2}\Delta_{\omega_{t_0}}^2(\tilde{\varphi}-\varphi)-\langle\varrho(\omega_t),dd^c(\tilde{\varphi}-\varphi)\rangle_t$, 
   and $w_{t_0,t}$ (uniformly) bounded at any order. 
   Uniform bounds at any order hold as well for the 
   $f_{\omega_t}^{\mathsf{Z}}=f_{\omega_X}^{\mathsf{Z}}+\mathsf{Z}\cdot\varphi_t
                             =f_{\omega}^{\mathsf{Z}}+t\mathsf{Z}\cdot(\tilde{\varphi}-\varphi)$;  
   these growth conditions near $D$ thus ensure us that 
    \begin{equation*}
     t\longmapsto \mathscr{F}_t:= \int_{\XD} \scal(\omega_t)f^{\mathsf{Z}}_{\omega_t}\,\frac{\omega_t^m}{m!}
    \end{equation*}
   is a smooth function of $t$, with derivative
    \begin{equation*}
     t\longmapsto \dot{\mathscr{F}}_t 
                  = \int_{\XD} \big(\dot{\scal}_t f^{\mathsf{Z}}_{\omega_t}
                    +           \scal(\omega_t) \big[\mathsf{Z}\cdot(\tilde{\varphi}-\varphi)\big]\big)\frac{\omega_t^m}{m!}
                    +\int_{\XD} \scal(\omega_t) f^{\mathsf{Z}}_{\omega_t}\,dd^c(\tilde{\varphi}-\varphi)\wedge\frac{\omega_t^{m-1}}{(m-1)!}, 
    \end{equation*}
   just as in the compact case. 
   And as in the compact case, 
   integrations by parts can be performed without boundary terms, 
   again thanks to the bounds mentioned above;  
   one thus ends with $\dot{\mathscr{F}}_t=0$ for all $t\in [0,1]$ (see e.g. \cite[Prop. 4.12.1]{gau}), 
   hence the result. 
  \cqfd
  
  \begin{rmk}
   The word ``character'' for the function $\mathscr{F}^D_{[\omega_X]}:\hD\to \R$ might appear slightly abusive, 
   as long as we have not checked that $\mathscr{F}^D_{[\omega_X]}([\mathsf{Z}_1,\mathsf{Z}_2])=0$ 
   for all $\mathsf{Z}_1,\mathsf{Z}_2\in \hD$. 
   As in the compact case, 
   this identity however follows at once from the invariance of $\mathscr{F}^D_{[\omega_X]}$ along 
   $\mathscr{M}^D_{[\omega_X]}$, and the stability of this class under automorphisms of $X$ parallel to $D$ and homotopic to $\id_X$. 
  \end{rmk}

\section{Link between smooth and Poincaré type Futaki characters}    \label{sctn_link}
 \subsection{Statement}                                              \label{subsec_sttmt}
  We keep the notations of the previous part; 
  in particular, $\omega_X$ is a smooth Kähler metric on $X$ compact, 
  and $\mathscr{F}_{[\omega_X]}^D:\mathfrak{h}_{\varparallel}^D \to \R$ denotes the Futaki character associated 
  to the space $\mathscr{M}_{[\omega_X]}^D$ of Poincaré type Kähler metrics on $\XD$ of class $[\omega_X]$. 
  
  Recall moreover that if $\mathsf{Z}\in \mathfrak{h}$, 
  we set $f^{\mathsf{Z}}_{\omega_X}$ for the normalised potential of its (Riemannian) gradient part, 
  relatively to $\omega_X$. 
  The purpose of this part is to state and prove an explicit relation between smooth and Poincaré type Futaki characters on $\hD$; 
  this is the main result of this note. 
  We use for this, as intermediates, 
  Futaki characters with respect to sub-divisors of $D$, 
  e.g. $D-D_j=\sum_{\ell=1,\ell\neq j}^N D_{\ell}$ if $D=\sum_{\ell=1}^N D_{\ell}$;  
  the Futaki character is denoted by $\mathscr{F}_{[\omega_X]}^{D-D_j}$ in this case, and is still defined on $\hD$.  
  We denote by $\mathscr{F}_{[\omega_X]}$ the usual Futaki character on $X$:
   \begin{prop}   \label{prop_futakis}
    For all $\mathsf{Z}\in \hD$ and for all $j=1,\dots, N$, 
    one has: 
     \begin{equation}    \label{eqn_futakis}
      \begin{aligned}
       \mathscr{F}_{[\omega_X]}^D(\mathsf{Z}) = \mathscr{F}_{[\omega_X]}^{D-D_j}(\mathsf{Z}) 
                                               +  4\pi \int_{D_j} f^{\mathsf{Z}}_{\omega_X} \frac{(\omega_X|_{D_j})^{m-1}}{(m-1)!}.                                                      
      \end{aligned}
     \end{equation}
    Consequently, for all $\mathsf{Z}\in \hD$, 
     \begin{equation}    \label{eqn_futakis1.5}
      \begin{aligned}
       \mathscr{F}_{[\omega_X]}^D(\mathsf{Z}) 
                  = \mathscr{F}_{[\omega_X]}(\mathsf{Z}) 
                    + 4\pi\sum_{j=1}^N \int_{D_j} f^{\mathsf{Z}}_{\omega_X} \frac{(\omega_X|_{D_j})^{m-1}}{(m-1)!}.                                                      
      \end{aligned}
     \end{equation}
   \end{prop}

 \subsection{Proof of Proposition \ref{prop_futakis}}               \label{subsec_PrfPropFutakis}
  Identity \eqref{eqn_futakis1.5} clearly follows from an inductive use of identity \eqref{eqn_futakis}, 
  the proof of which we focus on for the rest of this part. 
  
  Fix $Z\in \hD$. 
  To compute $\mathscr{F}^D_{[\omega_X]}(Z)$, 
  we first fix a Poincaré type Kähler metric $\omega\in \mathscr{M}^D_{[\omega_X]}$ as follows. 
  We take $\omega=\omega_X-dd^c\sum_{j=1}^N \log\big(-\log(|\sigma_j|^2_j)\big)$, 
  with $\sigma_j\in \mathscr{O}([D_j])$ such that $D_j=\{\sigma_j=0\}$, 
  and the $|\cdot|_j$ smooth hermitian metrics on the $[D_j]$, 
  chosen so that $\omega$ is indeed a (Poincaré type) metric on $\XD$ 
  -- see \cite[\S 1.1.1]{auv1} for details. 
  
  Fix now $j\in\{1,\dots,N\}$, set $\varphi_j=-\log\big(-\log(|\sigma_j|^2_j)\big)$, 
  $\psi_j=-\sum_{\ell\neq j}\log\big(-\log(|\sigma_j|^2_j)\big)$, 
  and define $\omega_t=\omega_X+dd^c(\psi_j+t \varphi_j)$ for $t\in[0,1]$. 
  Notice that 
  \textit{these are metrics of Poincaré type on $\XD$ for $t\in (0,1]$ only}, 
  as $\omega_{t=0}=\omega_X-dd^c\sum_{\ell\neq j}\log\big(-\log(|\sigma_j|^2_j)\big)$ 
  \textit{is of Poincaré type on $X\backslash(D-D_j)$} 
  -- assuming a good choice of the $|\cdot|_{\ell}$ for the positivity assertion. 
  Now by Proposition \ref{prop_InvceFutaki}, 
   \begin{equation}   \label{eqn_futaki_t}
    \mathscr{F}^D_{[\omega_X]}(Z) = \int_{\XD} \scal(\omega_t) f^{\mathsf{Z}}_{\omega_t}\, \frac{\omega_t^m}{m!}
   \end{equation}
  for all $t\in (0,1]$.   
  Observe however that the integrand tends uniformly to $\scal(\omega_0) f^{\mathsf{Z}}_{\omega_0}\, \frac{\omega_0^m}{m!}$ 
  away from $D_j$, 
  as $t$ goes to $0$. 
  Our strategy is hence to show that, 
  for the price of the correction
  $4\pi\int_{D_j} f_{\omega_X}^{\mathsf{Z}} \frac{(\omega_X|_{D_j})^{m-1}}{(m-1)!}$, 
  the formal limit $\int_{X\backslash (D-D_j)}\scal(\omega_0)f_{\omega_0}^{\mathsf{Z}}\,\frac{\omega_0^m}{m!}$ 
  is the limit of \eqref{eqn_futaki_t} as $t$ goes to 0; 
  in other words, we want to show that:
   \begin{equation}  \label{eqn_limit}
    \lim_{t\searrow 0 } \int_{\XD} \scal(\omega_t) f^{\mathsf{Z}}_{\omega_t}\, \frac{\omega_t^m}{m!} 
                        = \int_{X\backslash (D-D_j)}\scal(\omega_0)f_{\omega_0}^{\mathsf{Z}}\,\frac{\omega_0^m}{m!} 
                          +4\pi\int_{D_j} f_{\omega_X}^{\mathsf{Z}} \frac{(\omega_X|_{D_j})^{m-1}}{(m-1)!}, 
   \end{equation}
  which provides \eqref{eqn_futakis}, 
  by definition of $\mathscr{F}^{D}_{[\omega_X]}$ -- and its independence from $t\in(0,1]$ in \eqref{eqn_futaki_t} 
  --, 
  and of $\mathscr{F}^{D-D_j}_{[\omega_X]}$.   
  
  Set $D^j:=(D-D_j)|_{D_j}$; 
  admitting momentarily that 
   \begin{equation}    \label{eqn_equality1}
    \int_{D_j} f_{\omega_X}^{\mathsf{Z}} \frac{(\omega_X|_{D_j})^{m-1}}{(m-1)!}= 
    \int_{D_j\backslash D^j} f_{\omega_0}^{\mathsf{Z}} \frac{(\omega_0|_{D_j\backslash D^j})^{m-1}}{(m-1)!}, 
   \end{equation}
  our aim is to prove \eqref{eqn_limit} with 
  $\int_{D_j\backslash D^j} f_{\omega_0}^{\mathsf{Z}} \frac{(\omega_0|_{D_j\backslash D^j})^{m-1}}{(m-1)!}$ 
  instead of $\int_{D_j} f_{\omega_X}^{\mathsf{Z}} \frac{(\omega_X|_{D_j})^{m-1}}{(m-1)!}$. 
  The key point is the following technical lemma: 
   \begin{lem}   \label{lem_keylm}
    Let $f\in C^{\infty}\big(X\backslash(D-D_j)\big)$, 
    and $w\in C^{\infty}_{1}(X\backslash D_j)$. 
    Then 
     \begin{equation}   \label{eqn_keylm}
      \lim_{t\searrow 0 } \int_{\XD} \scal(\omega_t) (f+w)\, \frac{\omega_t^m}{m!} 
       = \int_{X\backslash(D-D_j)} \scal(\omega_0) (f+w)\, \frac{\omega_0^m}{m!} 
         + 4\pi\int_{D_j\backslash D^j} f \frac{(\omega_0|_{D_j\backslash D^j})^{m-1}}{(m-1)!}. 
     \end{equation}
    
   \end{lem}
  By ``$f\in C^{\infty}\big(X\backslash(D-D_j)\big)$'', we mean: 
  $f$ is smooth on $X\backslash(D-D_j)$, with derivatives bounded at any order with respect to any Poincaré type metric on $X\backslash(D-D_j)$, 
  e.g. $\omega_0$; 
  by ``$w\in C^{\infty}_{1}(X\backslash D_j)$'', 
  we mean: $w$ smooth on $X\backslash D_j$, 
  with derivatives at any order $\mathcal{O}\big(\big|\log|\sigma_j|_j\big|^{-1}\big)$
  with respect to any Poincaré type metric on $X\backslash D_j$, e.g. $\omega_X+dd^c\varphi_j$. 
  
  Lemma \ref{lem_keylm} is proven in next section. 
  Let us see for now how it applies to our situation. 
  One has: 
  $f^{\mathsf{Z}}_{\omega_t} = f^{\mathsf{Z}}_{\omega_0}+ t(\mathsf{Z}\cdot\varphi_j)$; 
  we already know that $f^{\mathsf{Z}}_{\omega_0}\in C^{\infty}\big(X\backslash(D-D_j)\big)$, 
  and we check easily that $(\mathsf{Z}\cdot\varphi_j)\in C^{\infty}_{-1}(X\backslash D_j)$ 
  thanks to the assumption that $\mathsf{Z}$ is parallel to $D_j$. 
  This way, by Lemma \ref{lem_keylm}, 
  $\int_{\XD} \scal(\omega_t) f_{\omega_0}^{\mathsf{Z}} \, \omega_t^m/m!$ tends to 
  $\int_{\XD} \scal(\omega_0) f_{\omega_0}^{\mathsf{Z}} \, \omega_0^m/m!\,
   +4\pi\int_{D_j\backslash D^j}f_{\omega_0}^{\mathsf{Z}} \, (\omega_0|_{D_j\backslash D^j})^{m-1}/(m-1)!$, 
  and 
  $\int_{\XD} \scal(\omega_t) (\mathsf{Z}\cdot\varphi_j) \, \omega_t^m /m!$ tends to 
  $\int_{\XD} \scal(\omega_0) (\mathsf{Z}\cdot\varphi_j) \, \omega_t^m /m!$ as $t$ goes to $0$ 
  -- all that matters here is actually this limit existing and being finite. 
  As a result, 
   \begin{align*}
    \int_{\XD} \scal(\omega_t) f_{\omega_t}^{\mathsf{Z}} \, \frac{\omega_t^m}{m!}
     &= \int_{\XD} \scal(\omega_t) f_{\omega_0}^{\mathsf{Z}} \, \frac{\omega_t^m}{m!}
       +t\int_{\XD} \scal(\omega_t) (\mathsf{Z}\cdot\varphi_j) \, \frac{\omega_t^m}{m!}                              \\
     &\xrightarrow{\,\, t \searrow 0 \,\,	} 
      \int_{X\backslash (D-D_j)} \scal(\omega_0) f_{\omega_0}^{\mathsf{Z}} \, \frac{\omega_0^m}{m!}
      +4\pi\int_{D_j\backslash D^j}f_{\omega_0}^{\mathsf{Z}} \, \frac{(\omega_0|_{D_j\backslash D^j})^{m-1}}{(m-1)!},
   \end{align*}
  as wanted. 
  
  Apart from the proof of Lemma \ref{lem_keylm}, 
  we are left with that of equality \eqref{eqn_equality1}. 
  We work on $D_j\backslash D^j$ 
  -- recall the notation $D^j=(D-D_j)|_{D_j}$ --, 
  where we set $\varpi_s=(1-s)(\omega_X|_{D_j})+s(\omega_0|_{D_j\backslash D^j})$; 
  these are Poincaré type metrics for $s>0$. 
  In the same fashion as in the proof of Proposition \ref{prop_HodgeDec}, 
  growths near $D^j$ allow us to say that 
   \begin{equation*}
    s \longmapsto \int_{D_j\backslash D^j} \big(f^{\mathsf{Z}}_{\omega_X}+s(\mathsf{Z}\cdot \psi_j)\big)\varpi_s^{m-1}
   \end{equation*}
  is smooth, 
  with derivative 
   \begin{equation*}
    \int_{D_j\backslash D^j} (\mathsf{Z}\cdot \psi_j)\,\varpi_s^{m-1}
     +(m-1) \int_{D_j\backslash D^j} \big(f^{\mathsf{Z}}_{\omega_X}+s(\mathsf{Z}\cdot \psi_j)\big)dd^c\psi_j\wedge\varpi_s^{m-2}.  
   \end{equation*}
  In order to conclude as in the proof of Proposition \ref{prop_HodgeDec}, 
  since $(\mathsf{Z}\cdot \psi_j)|_{D_j\backslash D^j}=(\mathsf{Z}|_{D_j\backslash D^j})\cdot(\psi_j|_{D_j\backslash D^j})$ 
  as $Z$ is parallel to $D_j$, 
  we check that the Hodge decomposition out of $D_j$ induces a Hodge decomposition on $D_j$, 
  up to the mean of the Riemannian/symplectic gradient potentials. 
  Namely, we check that 
   \begin{equation}   \label{eqn_HodgedecDj}
    \xi^{\mathsf{Z}|_{D_j}}_{\omega_X|_{D_j}}:=(\mathsf{Z}|_{D_j})^{\sharp_{(g_X|_{D_j})}}
                   =\xi_{\rm harm}|_{D_j}+d(f_{\omega_X}^{\mathsf{Z}}|_{D_j})+d^c(h_{\omega_X}^{\mathsf{Z}}|_{D_j}), 
   \end{equation}
  the extension to couples (Poincaré type metric $\varpi$ on $X\backslash(D-D_j)$, restriction of $\varpi$ on $D_j\backslash D^j$) 
  being dealt with as in Proposition \ref{prop_HodgeDec}. 
  Now, as harmonic 1-forms are exactly $d-$ and $d^c-$closed 1-forms on compact Kähler manifolds,  
  \eqref{eqn_HodgedecDj} is immediate from $\xi^{\mathsf{Z}|_{D_j}}_{\omega_X|_{D_j}}= \xi^{\mathsf{Z}}_{\omega_X}|_{D_j}$,  
  and this latter identity follows at once from the definition of $\mathsf{Z}|_{D_j}$. 
  Indeed, in local holomorphic coordinates $(z^1,\dots,z^m)$ such that $D_j$ is given by $z^1=0$, 
  write $\mathsf{Z}=\mathsf{Z}^k\frac{\partial}{\partial z^k}+\overline{\mathsf{Z}^k}\frac{\partial}{\partial \overline{z^k}}$, 
  and thus $\mathsf{Z}|_{D_j}=\mathsf{Z}^{\alpha}|_{D_j}\frac{\partial}{\partial z^{\alpha}}
                       +\overline{\mathsf{Z}}^{\alpha}|_{D_j}\frac{\partial}{\partial \overline{z^{\alpha}}}$ 
  -- we implicitly sum on repeated Latin indices over $\{1,\dots, m\}$, 
     and on Greek indices over $\{2,\dots,m\}$.  
  The dual 1-forms are given by: 
   \begin{equation*}
    \xi^{\mathsf{Z}}_{\omega_X} = \overline{\mathsf{Z}^{\ell}}(g_X)_{k\bar{\ell}}dz^k+\mathsf{Z}^{\ell}(g_X)_{\ell\bar{k}}d\overline{z^k}, 
    \,\, 
    \xi^{\mathsf{Z}|_{D_j}}_{\omega_X|_{D_j}} 
                                = \overline{\mathsf{Z}^{\beta}}|_{D_j}(g_X|_{D_j})_{\alpha\bar{\beta}}dz^{\alpha}
                                  +\mathsf{Z}^{\beta}|_{D_j}(g_X|_{D_j})_{\beta\bar{\alpha}}d\overline{z^{\alpha}}, 
   \end{equation*}
  hence the result after restriction to $D_j$ of $\xi^{\mathsf{Z}}_{\omega_X}$, 
  as $\mathsf{Z}^1|_{D_j}\equiv \overline{\mathsf{Z}^1}|_{D_j}\equiv 0$.  
  \cqfd

 \subsection{Main technical argument: proof of Lemma \ref{lem_keylm}}   \label{subsec_PrfKeyLm}
  \textit{Localisation of the problem. }
  Recall that $\omega_0=\omega_X+dd^c\psi_j$ is of Poincaré type on $X\backslash(D-D_j)$, 
  and that the $\omega_t=\omega_X+dd^c(t\varphi_j+\psi_j)$, $t\in(0,1]$, are of Poincaré type on $\XD$. 
  Now for all $t\in [0,1]$, 
  $\scal(\omega_t)\omega_t^m=2m\varrho(\omega_0)\wedge \omega_t^m
                             -m dd^c\log\big(\frac{\omega_t^m}{\omega_0^m}\big)\wedge\omega^{m-1}_t$. 
  On the one hand, 
  for $f$ and $w$ as in the statement, 
  as $(f+w)\varrho(\omega_0)\wedge \omega_t^m$ is uniformly dominated by $\omega^m$, 
   \begin{equation*}
    2m\int_{\XD}(f+w)\varrho(\omega_0)\wedge \omega_t^m 
      \to 2m\int_{\XD}(f+w)\varrho(\omega_0)\wedge \omega_0^m= \int_{\XD}\scal(\omega_0)(f+w)\,\omega_0^m
   \end{equation*}
  as $t$ tends to 0; 
  one recognises the first term in the right-hand side of \eqref{eqn_keylm}. 
  
  On the other hand,  
  thanks to the uniform convergence of 
  $dd^c\log\big(\frac{\omega_t^m}{\omega_0^m}\big)\wedge\omega^{m-1}_t$ 
  to $0$ far from $D_j$ (for $\omega_0$, say), 
  as $t$ tends to 0, 
  we can restrict to $f$ and $w$ with compact supports in  
  a neighbourhood $U$ of holomorphic coordinates $(z^1,\dots, z^m)$ centred at any point of $D_j$;  
  we also assume that $|z_{\ell}|\leq e^{-1}$ on $U$ for all $\ell$, 
  that $D_j\leq U=\{z^1=0\}$, 
  and that the possible other components of $D$ intersecting $U$ are respectively given by $\{z^2=0\},\dots,\{z^k=0\}$ 
  for the appropriate $k\in\{2,\dots, m\}$. 
  
  For fixed $t>0$, 
  we can write $\omega_t^m/\omega_0^m=v_t /[|z^1|^2\log^2(|z^1|^2)]$ on $U\backslash D$, 
  with $v_t$ positively bounded below, 
  and bounded up to order 2, for $\omega=\omega_{t=1}$; 
  these bounds are not uniform in $t$ though, 
  as $(\omega_t^m/\omega_0^m)\to 1$ far from $D_j$ when $t\searrow0$.  
  We rather write  
  $|\log(\omega_t^m/\omega_0^m)|\leq C+\log\big(1+t/[|z^1|^2\log^2(|z^1|^2)]\big)$ for a control uniform in $t$, 
  with $C>0$ independent of $t$. 
  
  Both controls come from the expansion 
  $\omega_t=\omega_0+t\frac{idz^1\wedge d\overline{z^1}}{|z^1|^2\log^2(|z^1|^2)}+\vareps_t$, 
  with $|\vareps_t|_{\omega}, |\nabla^{\omega}\vareps_t|_{\omega}, |(\nabla^{\omega})^2\vareps_t|_{\omega}\leq C t \big|\log|z^1|\big|^{-1}$, 
  where $C>0$ is independent of $t$. 
  
 ~
  
  \noindent
  \textit{Integration by parts. }
  Now as $dd^c\log(|z^1|^2)=0$ in $U\backslash D_j$, 
  for fixed $t>0$, 
   \begin{align*}
    \int_{U\backslash D} (f+w)dd^c\log\Big(\frac{\omega_t^m}{\omega_0^m}\Big)\wedge\omega^{m-1}_t
     = & \int_{U\backslash D} (f+w)dd^c\log\Big(\frac{v_t}{\log^2(|z^1|^2)}\Big)\wedge\omega^{m-1}_t  \\
     = & \int_{U\backslash D} \log\Big(\frac{v_t}{\log^2(|z^1|^2)}\Big)dd^c(f+w)\wedge\omega^{m-1}_t; 
   \end{align*}
  we perform this integration by parts without boundary terms, 
  as integrands are $L^1$ at every stage 
  (including the intermediate step, 
  where the integrand is 
  $d(f+w)\wedge d^c\log\big(\frac{v_t}{\log^2(|z^1|^2)}\big)\wedge\omega^{m-1}_t$). 
  
  Expand now $\omega_t^{m-1}$ as 
  $\omega^{m-1}+(m-1)t \frac{idz^1\wedge d\overline{z^1}}{|z^1|^2\log^2(|z^1|^2)}\wedge\omega_0^{m-2}+\tilde{\vareps}_t$, 
  with $|\tilde{\vareps}_t|_{\omega}\leq C t \big|\log|z^1|\big|^{-1}$; 
  this way, 
   \begin{equation}   \label{eqn_IntgrlSum}
    \begin{aligned}
     \int_{U\backslash D} &\log\Big(\frac{v_t}{\log^2(|z^1|^2)}\Big)dd^c(f+w)\wedge\omega^{m-1}_t  \\             
      = &\int_{U\backslash D} \log\Big(\frac{v_t}{\log^2(|z^1|^2)}\Big)dd^c(f+w)\wedge\omega_0^{m-1} \\
        &+(m-1)\int_{U\backslash D} \log\Big(\frac{v_t}{\log^2(|z^1|^2)}\Big)dd^c(f+w)
                                  \wedge\frac{t\,idz^1\wedge d\overline{z^1}}{|z^1|^2\log^2(|z^1|^2)} \wedge\omega_0^{m-2} \\
        & +\int_{U\backslash D} \log\Big(\frac{v_t}{\log^2(|z^1|^2)}\Big)dd^c(f+w)\wedge\tilde{\vareps}_t. 
    \end{aligned}
   \end{equation}
  We deal with the three summands of the right-hand side separately; 
  the aim is to show that when $t$ goes to 0, 
  the first summand provides the ``$\int_{D_j}$-term'' of \eqref{eqn_keylm}, 
  whereas the other two tend to 0. 
  
 ~ 
  
  \noindent
  \textit{First summand. }
  As $w|_{D_j}=0$,   
  (an easy adaptation of) the classical Lelong formula yields:
   $\int_{U\backslash D} \log(|z^1|^2)dd^c(f+w)\wedge\omega_0^{m-1}
       =-4\pi\int_{U\cap (D_j\backslash D^j)}f (\omega_0|_{D_j\backslash D^j})^{m-1}$. 
  Consequently, for $t>0$, as $\omega_t^m/\omega_0^m=v_t/[|z^1|^2\log^2(|z^1|^2)]$, 
   \begin{align*}
    \int_{U\backslash D} \log\Big(&\frac{v_t}{\log^2(|z^1|^2)}\Big)dd^c(f+w)\wedge\omega_0^{m-1} \\
     &= \int_{U\backslash D} \log\Big(\frac{\omega_t}{\omega_0}\Big)dd^c(f+w)\wedge\omega_0^{m-1}
        +4\pi\int_{U\cap (D_j\backslash D^j)}f (\omega_0|_{D_j\backslash D^j})^{m-1}.  
   \end{align*}
  The uniform controls 
  $|\log(\omega_t^m/\omega_0^m)|\leq C+\log\big(1+1/[|z^1|^2\log^2(|z^1|^2)]\big)$, 
  $|(dd^cf\wedge \omega_0^{m-1})/\omega_0^m |\leq C$, 
  $|(dd^cw\wedge \omega_0^{m-1})/\omega^m\big|\leq C\big|\log|z^1|\big|^{-1}$
  now allow us\footnote{the worst term to deal with is 
  $\int_{U\backslash D} \log\big(1+1/[|z^1|^2\log^2(|z^1|^2]\big)\big/\big|\log|z^1|\big|\,\omega^m$, 
  which is finite, 
  as $\log\big(1+ 1/[|z^1|^2\log^2(|z^1|^2]\big)\big/\big|\log|z^1|\big|=1+o(1)$ for $z^1$ small} 
  to argue by dominated convergence on the first summand of the right-hand side in the latter identity; 
  since the integrand tends to 0 as $t\searrow0$, 
  we get: 
   \begin{equation*}
    \lim_{t\searrow 0}\int_{U\backslash D}\, \log\Big(\frac{v_t}{\log^2(|z^1|^2)}\Big)dd^c(f+w)\wedge\omega_0^{m-1}
      = 4\pi\int_{U\cap (D_j\backslash D^j)}f (\omega_0|_{D_j\backslash D^j})^{m-1}.
   \end{equation*}
   
  \noindent
  \textit{Third summand of the right-hand side of \eqref{eqn_IntgrlSum}. } 
  Use the control on $\tilde{\vareps}_t$ to write:  
   \begin{align*}
    \bigg|\int_{U\backslash D} \log\Big(\frac{v_t}{\log^2(|z^1|^2)}&\Big)dd^c(f+w)\wedge\tilde{\vareps}_t\bigg|                 \\
     \leq& Ct\|dd^c(f+w)\|_{\omega} 
          \int_{U\backslash D}\Big|\log\Big(\frac{\omega_t^m}{\omega_0^m}\Big)+\log(|z^1|^2)\Big|\frac{\omega^m}{|\log(|z^1|^2)|};   
   \end{align*}
  the integral of the right-hand side is indeed finite 
  (same argument as in the footnote above), 
  and the left-hand side thus tends to 0 as $t\searrow0$. 
  
  ~
  
  \noindent
  \textit{Second summand of the right-hand side of \eqref{eqn_IntgrlSum}. } 
  This is probably the most delicate. 
  We rewrite the integral in play as 
   \begin{equation*}
    \int_{0<|z^1|\leq 1/e} \frac{t\,idz_1\wedge d\overline{z^1}}{|z^1|^2\log^2(|z^1|^2)}
     \int_{V_{z^1}} \log\Big(\frac{v_t}{\log^2(|z^1|^2)}\Big)\big(dd^c(f+w)\big)|_{V_{z^1}} \wedge(\omega_0|_{V_{z^1}})^{m-2}
   \end{equation*}
  where the $V_{z^1}$ are the slices $\{z^1={\rm constant}\}$ of $U\backslash D$. 
  On each such slice, 
  (the restriction of) $f+w$, $d(f+w)$ and $dd^c(f+w)$ are bounded, 
  with respect to (the restriction of) $\omega_0$, 
  hence $\int_{V_{z^1}}\big(dd^c(f+w)\big)|_{V_{z^1}} \wedge(\omega_0|_{V_{z^1}})^{m-2}=0$ for all $z^1\neq 0$. 
  Our integral can thus be rewritten as 
   \begin{equation*}
    \int_{0<|z^1|\leq 1/e} \frac{t\, idz_1\wedge d\overline{z^1}}{|z^1|^2\log^2(|z^1|^2)}
     \int_{V_{z^1}} \log(v_t)\big(dd^c(f+w)\big)|_{V_{z^1}} \wedge(\omega_0|_{V_{z^1}})^{m-2},
   \end{equation*}
  that is: 
   \begin{equation*}
    \int_{0<|z^1|\leq 1/e} idz_1\wedge d\overline{z^1}
     \int_{V_{z^1}} \frac{t\log[|z^1|^2\log^2(|z^1|^2)\cdot\omega_t^m /\omega_0^m]}{|z^1|^2\log^2(|z^1|^2)}\big(dd^c(f+w)\big)|_{V_{z^1}} 
                    \wedge(\omega_0|_{V_{z^1}})^{m-2}. 
   \end{equation*}
  Now for all $z^1\neq0$, $t\in(0,1]$, 
   \begin{align*}
     \Big|&\int_{V_{z^1}} \frac{t\log[|z^1|^2\log^2(|z^1|^2)\cdot\omega_t^m /\omega_0^m]}{|z^1|^2\log^2(|z^1|^2)}\big(dd^c(f+w)\big)|_{V_{z^1}} 
     \wedge(\omega_0|_{V_{z^1}})^{m-2}\Big|\\
      &\leq C\big\|\big(dd^c(f+w)\big)|_{V_{z^1}}\big\|_{\omega_0|_{V_{z^1}}}\vl(V_{z^1})
       \frac{t}{|z^1|^2\log^2(|z^1|^2)}\Big[1+\big|\log[t+|z^1|^2\log^2(|z^1|^2)]\big|\Big], 
   \end{align*}
  where $\vl(V_{z^1})=\int_{V_{z^1}}\omega_0^{m-1}$. 
  This volume, 
  as well as the supremums $\big\|\big(dd^c(f+w)\big)|_{V_{z^1}}\big\|_{\omega_0|_{V_{z^1}}}$ 
  are bounded below independently of $z^1$ (and of $t$!)
  -- notice that we restrict to directions parallel to $D_j$, 
  along which $\omega_0$ and $\omega$ are comparable. 
  Now, 
   \begin{align*}
   \int_{\{0<|z^1|\leq 1/e\}}&\frac{t\,idz^1\wedge d\overline{z^1}}{|z^1|^2\log^2(|z^1|^2)}
                          \Big[1+\big|\log[t+|z^1|^2\log^2(|z^1|^2)]\big|\Big] \\
     =&t \int_{\{0<|z^1|\leq 1/e\}}\frac{idz^1\wedge d\overline{z^1}}{|z^1|^2\log^2(|z^1|^2)}
      +t|\log t|  \int_{\{0<|z^1|\leq 1/e\}}\frac{idz^1\wedge d\overline{z^1}}{|z^1|^2\log^2(|z^1|^2)}\\
      & +\int_{\{0<|z^1|\leq 1/e\}}\frac{t}{|z^1|^2\log^2(|z^1|^2)}
                                   \log\Big(1+\frac{|z^1|^2\log^2(|z^1|^2)}{t}\Big)idz^1\wedge d\overline{z^1}. 
   \end{align*}
  As $t\searrow 0$, 
  the first two summands of the right-hand side clearly tend to 0; 
  as for the integrand of the third summand, 
  an elementary study of the function $x\mapsto x\log(1+1/x)$ on $(0,\infty)$ 
  shows that it is bounded above by 1, and tends to 0 as $t\searrow 0$. 
  A last use of dominated convergence thus gives that this third summand, 
  hence the whole second summand of \eqref{eqn_IntgrlSum}, tend to 0 as $t\searrow0$.  
  
 ~
  
  Summing up the above analysis of the three summands of the right-hand side of \eqref{eqn_IntgrlSum}, 
  we get: 
   \begin{equation*}
    \int_{U\backslash D} (f+w)dd^c\log\Big(\frac{\omega_t^m}{\omega_0^m}\Big)\wedge\omega_t^{m-1} 
     \xrightarrow{\,\, t\searrow 0\,\,} 4\pi  \int_{U\cap (D_j\backslash D^j)} f\, (\omega_0|_{(D_j\backslash D^j)})^{m-1}, 
   \end{equation*}
  and we saw this is equivalent to Lemma \ref{lem_keylm} for our (localised) $f$ and $w$. 
  \cqfd

\section{Application to extremal metrics of Poincaré type}                    \label{sctn_AppExtPK}
 \subsection{Extension of Proposition \ref{prop_futakis} (smooth divisor)}    \label{subsec_ExtPropFtks}
  Noticed that the integral term in \eqref{eqn_futakis} does not depend on the smooth metric $\omega_X\in \mathscr{M}_{[\omega_X]}$, 
  as neither $\mathscr{F}^{D}_{[\omega_X]}(\mathsf{Z})$ nor $\mathscr{F}^{D-D_j}_{[\omega_X]}(\mathsf{Z})$ do. 
  Considerations similar to those invoked when proving \eqref{eqn_equality1} tell us moreover 
  that for the price of replacing $D_j$ by $D_j\backslash D^j$, 
  one can replace $\omega_X$ by \textit{any} $\omega \in \mathscr{M}_{[\omega_X]}^{D-D_j}$,  
  $\omega|_{D_j\backslash D^j}$ being in that case an element of $\mathscr{M}_{[\omega_X]|_{D_j}}^{D^j}$. 
  
  One can go a step further, at least when the divisor is smooth, and take an $\omega \in \mathscr{M}_{[\omega_X]}^D$ 
  which is asymptotically a product near $D_j$, 
  i.e. for which there exist $a>0$, $\omega_j\in \mathscr{M}_{[\omega_X]|_{D_j}}$ and $\delta>0$ such that 
  as soon as $D_j=\{z^1=0\}$ in local holomorphic coordinates $(z^1,\dots, z^m)$, 
  then
   \begin{equation*}
    \omega = \frac{a\, idz^1\wedge d\overline{z^1}}{|z^1|^2\log^2(|z^1|^2)} + p^*\omega_j +\mathcal{O}\big(\big|\log|z^1|\big|^{-\delta}\big), 
   \end{equation*}
  where $p(z^1,\dots,z^m)=(z^2,\dots z^m)$, 
  and with the $\mathcal{O}$ understood at any order for $\omega$. 
  This way $\omega|_{D_j}$ still makes sense as an element of $\mathscr{M}_{[\omega_X]|_{D_j}}$, 
  as well as $f^{\mathsf{Z}}_{\omega}|_{D_j}$, and: 
   \begin{prop}[$D$ smooth]   \label{prop_futakis2}
    Let $\omega\in \mathscr{M}_{[\omega_X]}^D$, 
    and assume that $\omega$ is asymptotically a product near $D_j$, for $j\in\{1,\dots, N\}$. 
    Then for all $\mathsf{Z}\in \hD$, 
    one has: 
     \begin{equation}  
      \begin{aligned}
       \mathscr{F}_{[\omega_X]}^D(\mathsf{Z}) = \mathscr{F}_{[\omega_X]}^{D-D_j}(\mathsf{Z}) 
                                                 +  4\pi \int_{D_j} f^{\mathsf{Z}}_{\omega} \frac{(\omega|_{D_j})^{m-1}}{(m-1)!}.                                                      
      \end{aligned}
     \end{equation}
   \end{prop}
   
  \prf. --- 
   Assume that $\omega$ is asymptotically a product as above; 
   then $\omega=\omega_X+dd^c\big(\varphi+\tilde{\psi}\big)$, 
   with $\varphi=-a\log\big(-\log(|z^1|^2)\big)$, 
   and in local holomorphic coordinates $(z^1,\dots, z^m)$ such that $D_j=\{ z^1=0\}$, 
   $\tilde{\psi}=p^*\psi + \mathcal{O}\big(\big|\log|z^1|\big|^{-\delta}\big)$, 
   where the $\mathcal{O}$ is understood at any order for $\omega$, 
   and where $\psi\in C^{\infty}(D_j)$ is such that 
   $\omega_{D_j}^{\psi}:=\omega_X|_{D_j}+dd^c\psi\in \mathscr{M}_{[\omega_X]}$, 
   and $\omega|_{D_j}=\omega_{D_j}^{\psi}$.

   Taking $\mathsf{Z}\in \hD$,  
   $\mathsf{Z}\cdot\varphi=\mathcal{O}\big(\big|\log|z^1|\big|^{-1}\big)$ in coordinates as above, 
   so that $f^{\mathsf{Z}}_{\omega}= f^{\mathsf{Z}}_{\omega_{X}}+\mathsf{Z}\cdot\big(\varphi+\tilde{\psi}\big)$, 
   restricts to $f^{\mathsf{Z}}_{\omega}|_{D_j}+(\mathsf{Z}|_{D_j})\cdot \psi$ on $D_j$. 
   Now we know from the treatment of equality \eqref{eqn_equality1} in the proof of Proposition \ref{prop_futakis} 
   that $d\big(f^{\mathsf{Z}}_{\omega}|_{D_j}+(\mathsf{Z}|_{D_j})\cdot \psi\big)$ 
   is the gradient part in the Hodge decomposition of the dual 1-form of $(\mathsf{Z}|_{D_j})$ for $\omega_{D_j}^{\psi}$. 
   The analogue moreover holds when replacing $\psi$ by $t\psi$ for $t\in[0,1]$; 
   setting $\omega_{D_j}^t =\omega_X|_{D_j}+tdd^c\psi$ and $f_t=f^{\mathsf{Z}}_{\omega}|_{D_j}+t(\mathsf{Z}|_{D_j})\cdot \psi$, 
   we thus see that the derivative of $\int_{D_j} f_t\,(\omega_{D_j}^t)^{m-1}$ 
   vanishes thanks to the usual integration by parts, hence the result, 
   in view of \eqref{eqn_futakis}. 
    \cqfd

 \subsection{A numerical constraint on extremal metrics of Poincaré type}    \label{subsec_ExtCnstrnt}
  We apply what precedes to reformulate the numerical obstruction of \cite[\S4.2.2]{auv3}, 
  which is a constraint on \textit{extremal Poincaré type metrics} of class $[\omega_X]$: 
   \begin{thm}  \label{thm_ExtCnstt}
    Assume that there exists an extremal metric of Poincaré type of class $[\omega_X]$ on $\XD$, 
    and let $\mathsf{K}\in \hD$ be the Riemannian gradient of its scalar curvature. 
    Then for all $j=1,\dots, N$, setting $D^j=(D-D_j)|_{D_j}$, 
     \begin{equation}    \label{eqn_ExtCnstt}
      \sbar^D < \sbar_{D_j}^{D^j} +\frac{1}{4\pi \vl(D_j)}\big(\mathscr{F}_{[\omega_X]}^{D-D_j}(\mathsf{K})  - \mathscr{F}_{[\omega_X]}^D(\mathsf{K})\big),  
     \end{equation}
    where $\sbar^D$ (resp. $\sbar_{D_j}^{D^j}$) denotes the mean scalar curvature attached to $\mathscr{M}_{[\omega_X]}^D$ 
    (resp. to $\mathscr{M}_{[\omega_X]|_{D_j}}^{D^j}$). 
   \end{thm}
   
  \prf. --- Assume for a start that $D$ is smooth. 
   Let $\omega\in \mathscr{M}_{[\omega_X]}^D$ be extremal, and let $\mathsf{K}=\nabla \scal(\omega)\in \hD$, 
   where the (Riemannian) gradient $\nabla$ is computed with respect to (the Riemannian metric associated to) $\omega$. 
   According to \cite[Thm. 3]{auv3}, 
   $\omega$ is asymptotically a product near the divisor, 
   and induces an extremal metric $\omega_j \in \mathscr{M}_{[\omega_X]|_{D_j}}$ for all $j=1,\dots, N$. 
   We fix one of these $j$; 
   as $f^{\mathsf{K}}_{\omega}= \scal(\omega)-\sbar^{D}$, 
   Proposition \ref{prop_futakis2} implies: 
    \begin{equation*}    
      \begin{aligned}
       \mathscr{F}_{[\omega_X]}^D(\mathsf{K}) 
                         =& \mathscr{F}_{[\omega_X]}^{D-D_j}(\mathsf{K}) 
                                               +  4\pi \int_{D_j} (\scal(\omega)-\sbar^{D}) \frac{\omega_j^{m-1}}{(m-1)!}          \\
                         =& \mathscr{F}_{[\omega_X]}^{D-D_j}(\mathsf{K}) -4\pi\vl(D_j)\sbar^{D}
                                         +  4\pi \int_{D_j} \Big(\scal(\omega_j)-\frac{2}{a_j}\Big) \frac{\omega_j^{m-1}}{(m-1)!}  \\
                         =& \mathscr{F}_{[\omega_X]}^{D-D_j}(\mathsf{K}) -4\pi\vl(D_j)\Big(\sbar^{D}-\sbar_{D_j}+\frac{2}{a_j}\Big), 
      \end{aligned}
     \end{equation*}
   where $a_j \in (0,\infty)$ is such that: 
   given a neighbourhood of holomorphic coordinates $(z^1,\dots, z^m)$ in $X$ of any point of $D_j$ such that 
   $D_j$ locally corresponds to $z^1=0$, 
   then 
   $\omega= a_j\frac{i dz^1\wedge d\overline{z^1}}{|z^1|^2\log^2(|z^1|^2)}+p^*\omega_j+\mathcal{O}\big(\frac{1}{|\log(|z^1|)|^{\delta}}\big)$
   for some $\delta >0$, 
   and with $p(z^1,\dots, z^m)=(z^2,\dots,z^m)$. 
   As $a_j$ is positive, 
   one gets:  
    \begin{equation*}
     \mathscr{F}_{[\omega_X]}^{D-D_j}(\mathsf{K}) 
      >\mathscr{F}_{[\omega_X]}^D(\mathsf{K}) +4\pi\vl(D_j)(\sbar^{D}-\sbar_{D_j})  , 
    \end{equation*}
   of which \eqref{eqn_ExtCnstt} is simply a rewriting 
   -- as $D$ is smooth, $D^j=0$ on $D_j$.  
   
  ~

  \noindent
  \textit{The simple normal crossing case. }
   The asymptotically product behaviour of the extremal metric $\omega$ 
   is not clear anymore when the divisor admits (simple normal) crossings; 
   we thus content ourselves with applying Proposition \ref{prop_futakis}, 
   with $\mathsf{Z}=\mathsf{K}$, 
   and $\omega_X$ smooth, 
   and adapt our argument as follows. 
   Let $\varphi$ so that $\omega=\omega_X+dd^c\varphi$; 
   then $f^{\mathsf{K}}_{\omega}=f^{\mathsf{K}}_{\omega_X}+\mathsf{K}\cdot\varphi$, 
   that is, 
   $f^{\mathsf{K}}_{\omega_X}=f^{\mathsf{K}}_{\omega}-\mathsf{K}\cdot\varphi=\scal(\omega)-\sbar^D-\mathsf{K}\cdot\varphi$. 
   Remember that $f^{\mathsf{K}}_{\omega_X}$ is smooth on $X$, 
   and set for the following lines $\omega_{D_j}=\omega_X|_{D_j}$; 
   To compute $f^{\mathsf{K}}_{\omega_X}|_{D_j}$, 
   notice that by Remarks 4.4 and 4.7 in \cite{auv3}, 
   one can find ``tubes'' around $D_j$ in $\XD$ such that: $\scal(\omega)$ and $\mathsf{K}\cdot\varphi$ 
   tend uniformly on compact subsets of these tubes, respectively to $\scal(\omega_{D_j}+dd^c\psi)-2/a_j$ and $\mathsf{K}_{D_j}\cdot\psi$, 
   and where: $\psi$ is smooth on $D_j\backslash D^j$, 
   such that $\omega_{D_j}^{\psi}:=\omega_{D_j}+dd^c\psi\in \mathscr{M}_{[\omega_{D_j}]}^{D^j}$, 
   and $a_j>0$ is the inverse of the left-hand side of inequality (35) in \cite[Prop. 4.5]{auv3}.  
   
   As a consequence, 
   $f^{\mathsf{K}}_{\omega_X}|_{D_j\backslash D^j}=\scal(\omega_{D_j}+dd^c\psi)-\sbar^D-2/a_j-\mathsf{K}_{D_j}\cdot\psi$, 
   and Proposition \ref{prop_futakis} yields
    \begin{equation}  \label{eqn_sncFutakis}
     \mathscr{F}_{[\omega_X]}^D(\mathsf{K}) 
                         =  \mathscr{F}_{[\omega_X]}^{D-D_j}(\mathsf{K}) 
                             - 4\pi \int_{D_j\backslash D^j} \big(\scal(\omega_{D_j}^{\psi})-\sbar^D-\frac{2}{a_j}-\mathsf{K}_{D_j}\cdot\psi\big) 
                                               \frac{(\omega_{D_j})^{m-1}}{(m-1)!}.  
    \end{equation}
   We can be more specific when analysing $\omega$ near $D_j$, 
   and see that $\omega_{D_j}^{\psi}$ is extremal, 
   with $\mathsf{K}_{D_j}= \nabla\scal(\omega_{D_j}^{\psi})$ and $\nabla$ the Riemannian gradient with respect to $\omega_{D_j}^{\psi}$.  
   In other words, $f^{\mathsf{K}_{D_j}}_{\omega_{D_j}^{\psi}}=\scal(\omega_{D_j}^{\psi})-\sbar_{D_j}^{D^j}$, 
   hence 
    \begin{equation}   \label{eqn_fKDj}
     f^{\mathsf{K}_{D_j}}_{\omega_{D_j}}=\scal(\omega_{D_j}^{\psi})-\sbar_{D_j}^{D^j}-\mathsf{K}_{D_j}\cdot\psi. 
    \end{equation}
   As $\int_{D_j}f^{\mathsf{K}_{D_j}}_{\omega_{D_j}} \omega_{D_j}^{m-1}=0$ 
   by definition of the normalised holomorphic potential, 
   using \eqref{eqn_fKDj}, we can rewrite equation \eqref{eqn_sncFutakis} as: 
    \begin{equation*}
     \mathscr{F}_{[\omega_X]}^D(\mathsf{K}) 
                         =  \mathscr{F}_{[\omega_X]}^{D-D_j}(\mathsf{K}) 
                              - 4\pi\vl(D_j) \big(\sbar_{D_j}^{D^j}-  \sbar^D -\frac{2}{a_j} \big).     
    \end{equation*}
   We now conclude as in the smooth divisor case, using the positivity of $a_j$. 
   \cqfd

 \small

 \bibliographystyle{amsalpha}
 \bibliography{biblioPKFutaki}
   
 ~

\small \textsc{CMLA, École Normale Supérieure de Cachan, UMR 8536} \\
\indent 61 avenue du Président Wilson, 94230 Cachan, France        \\
\indent \url{hugues.auvray@ens-cachan.fr}

\end{document}